\documentclass[12pt]{article}
\usepackage{latexsym, amsfonts, amsmath, amsthm, amssymb,amscd}
\usepackage[dvips]{graphicx}
\usepackage{color}

\def\Fix{\mathop {\rm Fix}}
\def\OP{{\mathcal O}{\mathcal P}}
\def\F{{\mathcal F}}
\def\B{{\mathcal B}}
\def\Bm{ {\B}^*}
\def\M{{\mathcal M}}
\def\A{{\mathcal A}}
\def\ep{{\bf \epsilon}}
\def\mak{\mathop {mak}}

\def\I{\mathcal I}

\def\R{{\mathcal R}}
\def\C{{\mathcal C}}
\def\B{{\mathcal B}}
\def\W{{\mathcal W}}
\def\S{{\mathfrak S}}

\def\c{\mathbf{c}}
\def\r{\mathbf{r}}
\def\Rev{\mathop {\rm R}}
\def\DES{\mathop{DES}}

\def\ASC{\mathop{RISE}}
\def\des{\mathop{des}}

\def\maj{\mathop { maj}}
\def\mak{\mathop { mak}}
\def\stat{\mathop { stat}}

\def\bASC{\mathop {bRISE}}
\def\bDES{\mathop { bDES}}
\def\bDes{\mathop { bDes}}
\def\nbDes{\mathop { nbDes}}
\def\bMAJ{\mathop {bMAJ}}
\def\STAT{\mathop { STAT}}
\def\lsb{\mathop { lsb}}
\def\rsb{\mathop { rsb}}
\def\les{\mathop { les}}
\def\res{\mathop { res}}
\def\MAK{\mathop { MAK}}
\def\MIL{\mathop { MIL}}
\def\bmajMIL{\mathop { bmajMIL}}
\def\open{\operatorname{Open}}
\def\clos{\operatorname{Clos}}

\def\P{{\mathbb{P}}}
\def\Pm{{\P}^*}

\def\D{{\mathcal D}}
\def\Nn{{\mathbb{N}_n}}

\newtheorem{thm}{Theorem}[section]
\newtheorem{prop}[thm]{Proposition}

\newtheorem{cor}[thm]{Corollary}

\newtheorem{conj}[thm]{Conjecture}

\def\pf{\noindent {\it Proof.} }
\makeatletter
\newfont{\footsc}{cmcsc10 at 8truept}
\newfont{\footbf}{cmbx10 at 8truept}
\newfont{\footrm}{cmr10 at 10truept}
\makeatother 

\newtheorem{defn}[thm]{Definition}
\newtheorem{definition}[thm]{Definition}

\newtheorem{fact}[thm]{Fact}
\numberwithin{equation}{section}

\title{New Wilf-equivalence results for dashed patterns}
\date{}
\author{Anisse Kasraoui\thanks{The author was supported by grant no.\ 090038012 from the Icelandic
Research Fund and  grant S9607-N13 from Austrian Science Foundation FWF in the framework of the National Research Network
``Analytic Combinatorics and Probabilistic Number Theory{}''.}\\[0.5cm]
\centerline{\small Fakult\"at f\"ur Mathematik, Universit\"at Wien}\\ 
\centerline{\small Nordbergstra\ss e 15, A-1090 Vienna, Austria} \\
\centerline{\small\texttt{anisse.kasraoui@math.univ-lyon1.fr}}\\}

\begin{document}

\maketitle

\begin{abstract}
We give a sufficient condition  for the two dashed patterns 
$\tau^{(1)}-\tau^{(2)}-\cdots-\tau^{(\ell)}$ and $\tau^{(\ell)}-\tau^{(\ell-1)}-\cdots-\tau^{(1)}$  to be (strongly)
Wilf-equivalent. This permits to solve in a unified way several problems of Heubach and  Mansour 
on Wilf-equivalences on words and compositions, as well as a conjecture of Baxter and Pudwell on Wilf-equivalences on permutations.
We also give a better explanation of the equidistribution 
of the parameters $\MAK+\bMAJ$ and $\MAK'+\bMAJ$  on ordered set partitions.
These results can be viewed as consequences of a simple proposition  which states that the set valued statistics 
"descent set'' and "rise set'' 
are equidistributed over each equivalence class of the partially commutative monoid generated by a poset $(X,\leq)$.
\end{abstract}

%
%

\section{Introduction and Main results}

 The main purpose of this paper is to establish equidistribution properties about
dashed patterns in permutations, words and ordered set partitions. 
  
\subsection{Contribution to the Wilf-classification of dashed patterns}

Let $\P^*$ denote the free monoid generated by the set of positive integers $\P$. 
A \textit{dashed} (sometimes called  \textit{generalized} or \textit{vincular}) \textit{pattern} $p$ of length~$m$ is  
a word in $\Pm$ of length $m$ such that $p$ contains all letters in~$[\ell]:={1,2,\ldots,\ell}$ for a certain positive integer $\ell$ and 
in which two adjacent letters may or may not be separated by a dash. A dashed pattern $p=\tau^{(1)}-\tau^{(2)}-\cdots-\tau^{(\ell)}$ 
is said to be of type $(j_1, j_2,\cdots, j_{\ell})$ if the $\tau^{(i)}$ have lengths $j_i$. 
For example, $1\,3-2$ and $1\,1-3-4\,2$  are two dashed  patterns of type $(2,1)$ and $(2,1,2)$, respectively.  

 Let $w=w_1w_2\ldots w_n$ be a word in $\Pm$ and  
$p=p_1\ldots p_{j_1}-p_{j_1+1}\ldots p_{j_1+j_2}-\cdots-\ldots p_{m}$
 be a dashed pattern of length $m$. We say that the 
subsequence $w_{i_1}w_{i_2}\ldots w_{i_m}$, $1\leq i_1<i_2<\cdots<i_m\leq n$, is an occurrence 
of the pattern $p$ in~$w$ if \\[-0.8cm]
\begin{itemize}
 \item it is order isomorphic to $p$,  that is $w_{i_s}>w_{i_t}$ if and only $p_{s}>p_{t}$ for $1\leq s,t\leq m$, and \\[-0.8cm]
 \item the absence of a dash between two adjacent letters in the pattern $p$ indicates that
the corresponding letters in the subsequence $w_{i_1}w_{i_2}\ldots w_{i_m}$ must be adjacent in the word $w$.\\[-0.8cm]
\end{itemize}
For example, $w=2\,4\,1\,3\,5$ contains only two copies of $p=1-2\,3$, namely $2\,3\,5$ and $1\,3\,5$; although the subsequence $2\,4\,5$
of $w$ is isomorphic to $p$, it does not form a copy of $1-2\,3$ because 4 and 5 are not adjacent in $w$.

The main concern of this paper is the Wilf-equivalence on the set of dashed patterns.
 Let $\F=\bigsqcup_{n\geq 0} \F_n\subseteq \Pm$ be a collection of words such that $\F_n$, the set which consists of all 
words in $\F$ of length $n$, is finite for $n\geq 0$. The main collections we consider in this paper are:
\begin{itemize}
\item  The collection $\S=\bigsqcup_{n\geq 0}\S_n$ of permutations, where $\S_n$ is the symmetric group of order $n$;\\[-0.8cm]
\item  For $\ell$ a positive integer, the collection $\W_{[\ell]}=\bigsqcup_{n\geq 0}\W_{[\ell];n}$ of $\ell$-ary words (i.e., words the letters of which are $\leq \ell$), where $\W_{[\ell];n}$ stands for the set of words in $\W_{[\ell]}$ of length $n$;\\[-0.8cm]
\item  For a positive integer $s$  and a set $A\subseteq \P$, the collection $\C_{s,A}=\bigsqcup_{n\geq 0}\C_{s,A;n}$  
of integer compositions of $s$ all the parts of which are in $A$, where $\C_{s,A;n}$ is the set of compositions  in $\C_{s,A}$ of length $n$.
\end{itemize}
If $k$ is a nonnegative integer and $p$ is dashed pattern, 
we let $\F_n^{(k)}(p)$ denote the set of all
words in $\F_n$ containing exactly $k$ occurrences of the pattern $p$ and 
set $f_n^{(k)}(p)=|\F_n^{(k)}(p)|$. Then, two patterns $p$ and $q$ are said to be \textit{$\F$-Wilf-equivalent}, 
which is denoted $p \overset{\F}{\sim} q$, if $f_n^{(0)}(p)=f_n^{(0)}(q)$ for all $n\geq 0$. More generally, patterns $p$ and 
$q$ are said to be \textit{strongly $\F$-Wilf-equivalent} which is denoted $p\overset{\F}{\sim_{*}} q$ 
(or just $p\sim_{*} q$ if there is no confusion), 
if $f_n^{(k)}(p)=f_n^{(k)}(q)$ for all $n,k\geq 0$. 
The $\F$-Wilf-equivalence is obviously an equivalence relation and a natural problem is 
to determine its equivalence classes. While for classical patterns (i.e., dashed patterns all the 
letters of which are separated by a dash), 
much is known about $\S$-Wilf classes for patterns of length up to 7, 
for the dashed patterns, much less is known. Actually, we don't even know $\S$-Wilf classes 
for dashed patterns of length 4 (see e.g~\cite{St1}). The situation is similar for $\W_{[\ell]}$-Wilf equivalence 
and $\C_{s,A}$-Wilf equivalence. 

A related fundamental problem in the theory of patterns 
is to find necessary and sufficient conditions for two patterns to be Wilf-equivalent. 
In this paper, we give sufficient conditions for the two dashed patterns 
$\tau^{(1)}-\tau^{(2)}-\cdots-\tau^{(\ell)}$ and $\tau^{(\ell)}-\tau^{(\ell-1)}-\cdots-\tau^{(1)}$ to be Wilf-equivalent.
This leads to new Wilf-equivalences  on compositions, words and permutations, 
and permits to answer several problems listed by Heubach and Mansour about $\W_{[\ell]}$- and $\C_{s,A}$-Wilf equivalences 
as well as a conjecture of Baxter and Pudwell about $\S$-Wilf-equivalence.
  Before to present our result, we need some additional terminology.

 Recall that the \emph{descending runs} of a word $w=w_1w_2\ldots w_n\in\P^*$ are the maximal contiguous decreasing 
subsequences of $w$. For example, the word $w=3\,5\,4\,1\,6\,5\,5\,3\,6\,5$ 
have five descending runs: 3, 5\,4\,1, 6\,5, 5\,3 and 6\,5. 
Given a multiset $\M=\{D_1,D_2,\ldots,D_r\}$ of decreasing sequences, we let $\W(\M)$ denote
the set of words the descending runs of which are $D_1$, $D_2$,\ldots,$D_r$. 
For example, \\[-0.7cm]
\begin{itemize}
\item if $\M_1=\{3\,2\,1   \,,\,  6\,4  \,,\,  7\,5\}$, then 
$\W(\M_1)=\{ 3\,2\,1\,6\,4\,7\,5\,  \,,\,  3\,2\,1\,7\,5\,6\,4 \}$ but the permutation $6\,4\,3\,2\,1\,7\,5\notin\W(\I_1)$
since its descending runs are $6\,4\,3\,2\,1$ and $7\,5$;\\[-0.8cm]
\item if $\M_2=\{2\,1  \,,\,  2\,1  \,,\, 5\,3\}$, then $\W(\M_2)=\{ 2\,1\,2\,1\,5\,3\}$.
\end{itemize}
The \textit{reverse image} $\textbf{r}$ is the classical involutive transformation of $\Pm$ that maps each word
$w=x_1x_2\cdots x_n\in \Pm$ onto $\textbf{r}\,w= x_nx_{n-1}\cdots x_1$.
\begin{definition}
 A collection of words $\F\subseteq \Pm$ is said to be:
\begin{itemize} 
 \item \textit{run-complete} if for any multiset $\M$ of decreasing sequences, 
we have $\W(\M)\subseteq \F$ or $\W(\M)\cap \F=\emptyset$, i.e., 
if $\F$ contains a word in $\W(\M)$ then $\F$ contains every words in $\W(\M)$;\\[-0.8cm]
\item \textit{reverse-complete} if the reverse word $\r w$ is in $\F$ whenever $w \in \F$.
\end{itemize}
\end{definition}
It is easy to check that the collections of permutations $\S$, of $\ell$-ary words $\W_{[\ell]}$, of compositions $\C_{s,A}$,
are run-complete and reverse-complete.  In this paper, a special kind of dashed patterns 
will play an important role. Before to define them,
we need to introduce a partial order on $\Pm$.
\begin{definition}\label{def:partial order}
 Given two words $w$ and $w'$ in $\Pm$, $w$ is said to be \textit{below} $w'$
 if every letter in~$w$ is smaller than every letter in~$w'$. The partial order $\ll$
on $\Pm$ is defined as follows:  
\begin{align}
 w \ll w' \iff \text{$w=w'$ or $w$ is below $w'$.}
\end{align}
Two words $w$ and $w'$ are said to be \textit{comparable} if  $w\ll w'$ or $w\gg w'$, 
otherwise they are \textit{incomparable}.
\end{definition}

\begin{definition}
 A dashed pattern $p=\tau^{(1)}-\tau^{(2)}-\cdots-\tau^{(\ell)}$  is said to be\\[-0.7cm] 
\begin{itemize}
\item[-] \emph{connected} if, for $i=1,2,\ldots,\ell-1$, $\tau^{(i)}$ and $\tau^{(i+1)}$ are incomparable or equal;\\[-0.8cm] 
\item[-] \emph{piecewise decreasing} (resp., \textit{increasing}) if each $\tau^{(i)}$ is  (strictly)  decreasing (resp., increasing).
\end{itemize}
\end{definition}
For example, the  pattern $p=2\,3\,4-1-1\,2\,4$  is piecewise increasing but not connected since the two first 
sequences  are comparable ($2\,3\,4\gg 1$). The pattern  $q=5\,2-1\,4-3$ is connected but not piecewise  increasing (resp., decreasing).
The pattern  $q=5\,2-4\,1-3$ is piecewise decreasing and connected. 

If $p_1, p_2, \ldots, p_m$ are dashed patterns, the $m$-statistic which associates to a word $w\in\Pm$ 
the vector $(a_1,a_2,\cdots,a_m)$ where $a_i$ is the number of occurences of  $p_i$ in $w$ will be denoted 
by $(p_1,\ldots, p_m)$. We say that two vectors $(p_1, p_2, \ldots, p_m)$ and $(q_1, q_2, \ldots, q_m)$ of
dashed patterns are strongly $\F$-Wilf-equivalent, which is denoted 
$(p_1, p_2, \ldots, p_m)\overset{\F}{\sim_{*}} (q_1, q_2, \ldots, q_m)$, 
 if and only if the $m$-statistics $(p_1,\ldots, p_m)$ and $(q_1, \ldots, q_m)$ 
have the same distribution over~$\F_n$ for all integers $n\geq0$, that is 
for all integers $n,k_1,k_2,\ldots,k_m\geq 0$,
we have
$$\big(f_n^{(k_1)}(p_1),f_n^{(k_2)}(p_2),\ldots,f_n^{(k_m)}(p_m)\big)
=\big(f_n^{(k_1)}(q_1),f_n^{(k_2)}(q_2),\ldots,f_n^{(k_m)}(q_m)\big).$$
Consider the two transformations $\Rev$ and $<\r|\,\cdot>$ of dashed patterns defined for 
 $p=\tau^{(1)}-\tau^{(2)}-\cdots-\tau^{(\ell)}$ by 
\begin{align}
 \Rev p&:=\tau^{(\ell)}-\tau^{(\ell-1)}-\cdots-\tau^{(1)},\\
 <\r|\,p>&:=\r \tau^{(1)}-\r \tau^{(2)}-\cdots-\r \tau^{(\ell)}.
\end{align}
For example, if $p=2\,5\,4-3-1\,2$, we have $\Rev p=1\,2-3-2\,5\,4$ and $<\r|\,p>=4\,5\,2-3-2\,1$.
The following is the main result of this section.
\begin{thm}\label{thm:Wilf-equiv}
Let $\F\subseteq \Pm$ be a run-closed collection of words.  
Then, for any piecewise decreasing and connected patterns $p_1, p_2,\ldots,p_m$, we have
\begin{align}
 (p_1,p_2,\ldots,p_m) \overset{\F}{\sim_{*}} (\Rev p_1,\Rev p_2,\ldots,\Rev p_m).
\end{align}
In particular, the joint distribution of $(p_i,\Rev p_i)$ over $\F_n$ is symmetric for any integer $n\geq0$.
If, in addition, $\F$ is reverse-closed, then
\begin{align}
 (p_1,p_2,\ldots,p_m) \overset{\F}{\sim_{*}} (<\r|\,p_1>,<\r|\,p_2>,\ldots,<\r|\,p_m>).
\end{align}  
\end{thm}

The above result permits to unify previous results on Wilf-equivalences as well as to solve several problems on Wilf-equivalences.
For instance, since the collection of compositions $\C_{s,A}$
is run- and reverse-complete, we obtain the following result  
as an immediate consequence of Theorem~\ref{thm:Wilf-equiv}.

\begin{cor}\label{cor:Wilf-equivalences C}
 For any set $A\subseteq P$ and any positive integer $s$, 
we have the following (strong) $\C_{s,A}$-Wilf-equivalences:\\[0.3cm]
(1) type $(2,1)$ patterns:
\begin{align*}
(a)\;\; 1\,2-2\sim_{*} 2\,1-2\qquad (b)\;\; 1\,3-2\sim_{*} 3\,1-2.
\end{align*}
(2) type $(3,1)$ patterns:
\begin{align*}
&(a)\;\; 1\,2\,3-1\sim_{*} 3\,2\,1-1\qquad (b)\;\; 1\,2\,3-2\sim_{*} 3\,2\,1-2 \qquad (c)\;\; 1\,2\,3-3\sim_{*} 3\,2\,1-3\\
&(d)\;\; 1\,2\,4-3\sim_{*} 4\,2\,1-3\qquad (e)\;\; 1\,3\,4-2\sim_{*} 4\,3\,1-2. 
\end{align*}
(3) type $(2,2)$ patterns:
\begin{align*}
&(a)\;\; 1\,3-1\,2\sim_{*} 1\,2-1\,3\qquad (b)\;\;  1\,2-2\,3\sim_{*} 2\,1-3\,2 \qquad (c)\;\; 1\,3-2\,3\sim_{*} 2\,3-1\,3\\
&(d)\;\; 1\,4-2\,3\sim_{*} 2\,3-1\,4\qquad (e)\;\; 1\,3-2\,4\sim_{*} 2\,4-1\,3. 
\end{align*}
(4) type $(1,1,2)$ patterns:
\begin{align*}
&(a)\;\; 1-1-1\,2\sim_{*} 1-1-2\,1\qquad (b)\;\;  2-2-1\,2\sim_{*} 2-2-2\,1\\
&(c)\;\; 2-2-1\,3\sim_{*} 2-2-3\,1. 
\end{align*}
(5) type $(1,2,1)$ patterns:
\begin{align*}
&(a)\;\; 1-1\,2-2\sim_{*} 1-2\,1-2\qquad (b)\;\; 1-1\,3-2\sim_{*} 1-3\,1-2\\
&(c)\;\; 2-1\,3-3\sim_{*} 2-3\,1-3\qquad (d)\;\; 2-1\,4-3\sim_{*} 2-4\,1-3. 
\end{align*}
\end{cor}
It is worth noting that the above $\C_{s,A}$-Wilf-equivalences for $(2,1)$ patterns were first obtained by Heubach 
and Mansour (see Theorem~{5.38} and Theorem~{5.19} in~\cite{HeMa}) and their proof is far to be obvious. All the other
$\C_{s,A}$-Wilf-equivalences in Corollary~\ref{cor:Wilf-equivalences C} settle problems 
of Heubach and Mansour (see Questions (1-a), (2), (3) and (4) in~\cite[pp. 177--178]{HeMa}).

Since the collections of $\ell$-ary words $\W_{[\ell]}$ and of permutations $\S$ are run- and reverse-complete,
all the Wilf-equivalences given in Corollary~\ref{cor:Wilf-equivalences C} are also $\W_{[\ell]}$-Wilf-equivalences  
and $\S$-Wilf-equivalences. However, several of these equivalences are trivial.
Recall that the symmetry class of a dashed pattern $p$ is the set of patterns $\{p,\r\, p,\c\, p,\r\,\c\,p\}$,  
where $\r$ is the reverse map (that we have defined previously) and $\c$ is the complement transformation,
defined for a word $w=w_1w_2\ldots w_n$ by $\c\,w=y_1 y_{2}\ldots y_n$ where $y_i=M+1-w_i$, with $M$ being the maximum
of the letters of $w$.
For example, the patterns $2-31$, $2 -13$, $13-2$
and $31-2$ form a  symmetry class. It is easy to see (and well-known) that two dashed patterns 
in the same symmetry class are (strongly) $\S$- and $\W_{[\ell]}$-Wilf equivalent.
We summarize the non-trivial $\S$- and $\W_{[\ell]}$-Wilf equivalences of length 4 which can be deduced
from Theorem~\ref{thm:Wilf-equiv} (or Corollary~\ref{cor:Wilf-equivalences C}) in the two following results.

\begin{cor}\label{cor:Wilf-equivalences W}
 For any positive integer $\ell$, 
we have the following (strong) $\W_{[\ell]}$-Wilf-equivalences:
\begin{align*}
&(\rm{i})\;\; 1\,3-1\,2\sim_{*} 1\,2-1\,3\qquad (\rm{ii})\;\;  1\,2-2\,3\sim_{*} 2\,1-3\,2 \qquad 
(\rm{iii})\;\; 1\,3-2\,4\sim_{*} 2\,4-1\,3\\
&(\rm{iv})\;\; 1\,2\,3-1\sim_{*} 3\,2\,1-1 \;(\sim_{*} 1\,2\,3-3) \quad 
(\rm{v})\;\; 1\,2\,4-3\sim_{*} 4\,2\,1-3 \;(\sim_{*} 1\,3\,4-2)\\
&(\rm{vi})\;\; 1-1\,2-2\sim_{*} 1-2\,1-2\qquad \qquad\;(\rm{vii})\;\; 1-1\,3-2\sim_{*} 1-3\,1-2\\
&(\rm{viii})\;\; 2-1\,4-3\sim_{*} 2-4\,1-3\qquad\qquad (\rm{ix})\;\; 1\,2-1-1\sim_{*} 2\,1-1-1 \;(\sim_{*} 1\,2-2-2),
\end{align*}
where parentheses indicate trivial equivalences (obtained from symmetry classes).
\end{cor}
 All the $\W_{[\ell]}$-Wilf-equivalences in Corollary~\ref{cor:Wilf-equivalences W} settle problems 
of Heubach and Mansour (see Questions (2)--(5) in~\cite[pp. 238--239]{HeMa}). 
We also have 3 new (strong) $\S$-Wilf-equivalences for patterns of length 4. 
\begin{cor}\label{cor:Wilf-equivalences S}
 We have the following (strong) $\S$-Wilf-equivalences:
\begin{align*}
(a)\;\; 1\,2\,4-3\sim_{*} 4\,2\,1-3\qquad (b)\;\; 2-1\,4-3\sim_{*} 2-4\,1-3\qquad (c)\;\; 1\,3-2\,4\sim_{*} 2\,4-1\,3.
\end{align*}
\end{cor}
It is worth noting that the $\S$-Wilf-equivalence $(a)$ in Corollary~\ref{cor:Wilf-equivalences S} settles Conjecture~17~(b) in~\cite{BaPu}.
Baxter also obtained a proof which will be published in a upcoming paper (see the remarks after Conjecture~17 in~\cite{BaPu}).

Before to conclude this section, we want to notice that Theorem~\ref{thm:Wilf-equiv} can be obtained 
from Clarke, Steingr\'{\i}msson and Zeng's transformation $\ep$ on words~\cite[Section 6]{ClStZe}, 
which can be seen as a particular case of the transformation $\theta$ of Foata and Randrianarivony~\cite{FoRa}. 
It was shown in~\cite{ClStZe} that this transformation $\ep$ has two interesting properties: it preserves the 
descending runs and it exchanges the combinatorial parameters $\les$ and $\res$ 
(it is easy to see that we have $\les=(3\,1-2)$ and $\res=(2-3\,1)$). Actually, it is not difficult to see 
that the transformation~$\ep$ also exchanges the parameters $(p)$ and $(\Rev p)$ for any piecewise decreasing and connected pattern $p$
(see Theorem~\ref{thm:epsilon}). 
This leads to the following refinement of Theorem~\ref{thm:Wilf-equiv}.
\begin{thm}\label{thm:Wilf-equiv-local}
For any finite multiset $\M\subseteq \B$ and 
any piecewise decreasing and connected patterns $p_1, p_2,\ldots,p_m$, 
the $m$-tuple $(p_1,p_2,\ldots,p_m)$ and $(\Rev p_1,\Rev p_2,\ldots,\Rev p_m)$
have the same distribution over $\W(\M)$.
\end{thm}

We will show in this paper that Theorem~\ref{thm:Wilf-equiv-local} is a particular case 
of a result on partially commutative monoid which also permits to obtain 
a more simple and deeper explanation of an equidistribution 
result on ordered set partitions.  
\subsection{Statistics on ordered set partitions}
 Recall that an \textit{ordered set partition} $\pi=(B_1,B_2,\cdots ,B_k)$ 
of $\Nn:=\{1,2,\ldots,n\}$ is a sequence  of disjoint and nonempty subsets $B_i$, called blocks, whose union is $\Nn$. 
The set of all ordered set partitions of $\Nn$ into $k$ blocks is denoted $\OP_{n}^k$. By convention,
a block of a partition will be represented by the decreasing rearrangement of its elements and the blocks of a partition
will be separated by a vertical line.
For instance, the partition $\big(\{1\},\{5,8\},\{2,6,9\},\{3\},\{4\,7\}\big)\in \OP_9^5$ will be written 
as $1\,|\,8\,5\,|\,9\,6\,2\,|\,3\,|\,7\,4$.
It is well-known that $|\OP_{n}^k|=k!\,S(n,k)$ where $S(n,k)$ is the $(n,k)$-th Stirling number of the second kind
(see e.g.~\cite{Aigner}). The 
natural $q$-analogue $[k]_q!S_q(n,k)$ of the counting function $k!\,S(n,k)$ have arised as the generating
function for the distribution of several statistics on ordered set partitions (see e.g.~\cite{St,IsKaZe,KaZe}). 
Here, $[k]_q!$ is  the usual $q$-factorial  $[1]_q[2]_q\cdots[k]_q$, 
where $[j]_q=1+q+q^2+\cdots+q^{j-1}$, and $S_q(n,k)$ is the
$q$-Stirling of the second kind  defined by 
\begin{align*}
 S_q(n,k)=q^{k-1}\,S_q(n-1,k-1)+[k]_q\,S_q(n-1,k)
\end{align*}
for $n\geq k\geq 1$, and $S_q(n,k)=\delta_{n,k}$ if $n$ or $k$ are $0$.
The systematic study of statistics on ordered set partitions has its origins in the work of Steingr\'{\i}msson~\cite{St}.
Following Steingr\'{\i}msson, a statistic $\stat$ on ordered set partitions  such that
\begin{align*}
 \sum_{\pi\in\OP_n^k}q^{\stat \,\pi}= [k]_q!S_q(n,k) \quad(n\geq k\geq 1),
\end{align*}
is said to be \textit{Euler-Mahonian}.
Steingr\'{\i}msson~\cite{St} found several Euler-Mahonian statistics on ordered set partitions, most of them
reflect quite naturally the recursion for the $q$-Stirling numbers, but some others inspired by permutations statistics
don't have  simple recursive structures; they are qualified as ``\textit{hard}''. We are mainly concerned in this paper 
with the two hard statistics $\MAK+\bMAJ$ and ${\MAK}'+\bMAJ$. They are defined as follows.
 Suppose we are given a partition $\pi=B_1|B_2|\cdots |B_k\in\OP_n^k$.\\
-- An integer $i$, $1\leq i\leq k-1$, is said to be a \emph{block descent} (resp., \textit{block rise})
if $B_i \gg B_{i+1}$ (resp., $B_i \ll B_{i+1}$), where  $\gg$ is the partial order introduced in
Definition~\ref{def:partial order}. The sets of block descents and rises of $\pi$ will be 
denoted $\bDES(\pi)$ and $\bASC(\pi)$. The \emph{block major index} of $\pi$, denoted by $\bMAJ \pi$,
is defined as the sum of the block descents in $\pi$. \\
-- The \textit{opener} of a block is its least element 
and the \textit{closer} is its greatest element. We will denote by $\open(\pi)$ and
$\clos(\pi)$ the sets of openers and closers of the blocks of~$\pi$, respectively. \\
-- For an integer $i$, $1\leq i\leq n$, we let $\rsb_i \,\pi$ (resp., $\lsb_i\, \pi$) denote the number of blocks~$B$ 
in $\pi$ to the right (resp., left) of the block containing $i$ such that the opener of $B$ is smaller than~$i$ 
and the closer of $B$ is greater than $i$.
We then define the statistics $\rsb$ and $\lsb$ as the sum of their
coordinate statistics, i.e.
$$
\rsb \pi=\sum_{i=1}^n{\rsb}_{i}\,\pi \quad\text{and}\quad \lsb \pi=\sum_{i=1}^n{\lsb}_{i}\,\pi.
$$ 
-- The partition statistics $\MAK$ and ${\MAK}'$ are then 
defined by 
\begin{align*}
\MAK \pi&= \rsb\,\pi+\sum_{i\in\, \clos(\pi)}(n-i)
\quad\text{and}\quad
{\MAK}' \pi= \rsb\,\pi+\sum_{i\in\, \open(\pi)}(i-1).
\end{align*}
For example, if $\pi=8\,5\,|\,1\,|\,9\,6\,2\,|\,7\,4\,|\,3$, we have\\ 
-- $\bDES(\pi)=\{1,4\}$ and ${\bASC(\pi)=\{2\}}$, whence  $\bMAJ \pi=1+4=5$;\\
-- $\open(\pi)=\{1,2,3,4,5\}$ and $\clos(\pi)=\{1,3,7,8,9\}$;\\
-- $\big(\rsb_i\, \pi\big)_{1\leq i\leq n}=(0,0,0,0,2,1,0,1,0)$ and 
$\big(\lsb_i\, \pi\big)_{1\leq i\leq n}=(0,0,1,1,0,1,2,0,0)$, whence ${\rsb}\,\pi=4$ and ${\lsb}\,\pi=5$;\\
-- $\MAK \pi= 5+(8+6+2+1+0)=19$ and ${\MAK}' \pi= 5+(0+1+2+3+4)=15$.\\ [-0.3cm]

\textbf{Result A}. \textit{For $n\geq k\geq 1$, the partition statistics 
$\MAK+\bMAJ$ and ${\MAK}'+\bMAJ$ are equidistributed over~$\OP_n^k$}.\\[-0.3cm]

Result A, originally conjectured by Steingr\'{\i}msson, was proved by Zeng and the author
(see Theorem~3.3 in~\cite{KaZe}). An interesting fact is that $\MAK$ and ${\MAK}'$ can be seen as extension of 
the Mahonian permutation statistics $\mak$ (introduced by Foata and Zeilberger~\cite{FoZe}) and its 
variant ${\mak}'$ (introduced by Clarke, Steingr\'{\i}msson and Zeng~\cite{ClStZe}). More precisely, if $\sigma\in\S_n$, let 
$D_1,D_2,\ldots,D_k$ be the descending runs  of $\sigma$ (listed from left to right) and recall that 
$\des(\sigma)$ is the \textit{number of descents} of $\sigma$, where a descent is as usual an integer $i$ such 
that $\sigma(i)>\sigma(i+1)$.  Then the statistics 
$\mak$ and $\mak '$ can be defined by
\begin{align}
 \mak \sigma& :=\MAK (D_1\,|\,D_2\,|\cdots|\,D_k)+{n+1\choose 2}-kn\label{eq:MAK-mak}\\
{\mak} ' \sigma& :={\MAK}' (D_1\,|\,D_2\,|\cdots|\,D_k)+{n+1\choose 2}-kn.\label{eq:MAK'-mak'}
\end{align}
For example, we have $\mak(1\,8\,5\,9\,6\,2\,3\,7\,4)=\MAK(1\,|\,8\,5\,|\,9\,6\,2\,|\,3\,|\,7\,4)=19$. \\[-0.3cm]

\textbf{Result B}. \textit{For $n\geq 1$, the permutation statistics $(\des,\mak)$ and $(\des,\mak')$ 
are equidistributed over~$\S_n$}.\\[-0.3cm]

 Result B was obtained by Clarke, Steingr\'{\i}msson and Zeng (see Proposition~16 in~\cite{ClStZe})
in their study of Euler-Mahonian statistics on permutations and its proof relies on a simple
transformation on $\S_n$. This contrasts with the proof of Result A in~\cite{KaZe} which is based on a non trivial path model 
for ordered set partitions (see Sections~8 and~9 in~\cite{KaZe}) and showed no connection with Result~B. 
Altogether, this leads to two natural questions:

-- Since $\MAK$ and ${\MAK}'$ are extensions of $\mak$ and $\mak'$, can Result~A and Result~B
be unified?

-- Is there a simple proof of Result B?\\[-0.3cm]
 
We will answer these questions by giving a simple proof of a refinement of Result~B
which generalizes Result A.   

\begin{thm}\label{thm:application1}
 The 3-tuple of parameters 
$$(\MAK,{\MAK}',\bDES)\quad\text{and}\quad({\MAK}',\MAK,\bDES)$$
are equidistributed over $\OP_n^k$. 
\end{thm}
Notice that in order to see that Theorem~\ref{thm:application1} implies Result~A, 
it suffices to observe that via its descending runs, a permutation can be seen 
as an ordered set partition with no block descent.  Before to conclude this section, we want to indicate to the reader that the ``connections'' between 
Euler-Mahonian partition and permutation statistics is far to be understood. In the last section of this paper, 
we present a conjectured equidistribution result on ordered set partitions that generalizes an important 
equidistribution property over the set of permutations.
 
\subsection{The key result}

  Although it seems a priori that the two topics we discuss are disconnected, they are 
in fact strongly related if we work in the partially commutative monoid generated by 
the poset $\big(\mathcal{P}_f(\P),\ll \big)$, where $\mathcal{P}_f(\P)$ is the collection of all 
finite subsets of $\P$ and $\ll$  is the partial order introduced in Definition~\ref{def:partial order}. 
All the results presented in this paper can be derived from the following result on  
partially commutative monoid.
\begin{thm}\label{thm:EquiAscDesA}
For any poset $(X,\leq)$,  the  set-valued statistics ``\textit{set of descents}'' and ``\textit{set  of rises}''
are equidistributed over each equivalence class of the partially commutative monoid $L(X,\leq)$.
\end{thm}
In Section~2, we recall some definitions in the theory of partially commutative monoid 
and we see how Theorem~\ref{thm:EquiAscDesA} implies all the results presented in this paper.
In Section~3, we prove Theorem~\ref{thm:EquiAscDesA}. 
We end this paper with some remarks and problems.

%
%
\section{The partially commutative monoid $L(\B,\ll)$}


\subsection{Partially commutative monoids}

 We first recall the construction of the partially commutative monoid $L(X,\leq)$
generated by a nonempty poset~$(X,\leq)$. Let  $X^*$ be the free monoid generated by $X$. 
Two words~$w$ and~$w'$ in $X^*$ are said to be \textit{adjacent}  if there exist two words $u$ and $v$ 
and an ordered pair~${(a,a')\in X^2}$ of distinct and comparable elements, i.e $a <a'$ or $a > a'$,
such that~${w=uaa'v}$ and~$w'=ua'av$. They are said to be \textit{equivalent} if they are equal,
or if there exists a sequence of words~$w_0$,$w_1$,\ldots,$w_p$
such that~$w_0=w$,~$w_p=w'$ and~$w_{i-1}$ and~$w_{i}$ are adjacent for $1\leq i\leq p$. This defines an equivalence
relation $R_{\leq}$ on $X^*$, compatible with the multiplication in~$X^*$. 
Then, $L(X,\leq)$ is defined as the quotient monoid $X^*/R_{\leq}$. The equivalence class 
of a word $w\in X^*$ will be denoted by $\big[w\big]$.

 We now recall the definition of \textit{descent} and \textit{rise} (sometimes called \textit{ascent}) in  a word. 
Let  $w=x_1x_2\cdots x_n$ be a word of length $n$ in $X^*$. Then, the integer $i$, $1\leq i\leq n-1$, is said to 
be a descent (resp., rise) if $x_i>x_{i+1}$ (resp., $x_i<x_{i+1}$). The \textit{descent set}~$\DES(w)$ 
and \textit{rise set} $\ASC(w)$ of $w$ are
\begin{align*}
 \DES(w)&=\{i\,:\,1\leq i\leq n-1\quad\text{and}\quad x_i>x_{i+1}\}\\
 \ASC(w)&=\{i\,:\,1\leq i\leq n-1\quad\text{and}\quad x_i<x_{i+1}\}.
\end{align*}
A word $w=x_1x_2\cdots x_n\in X^*$ is said to be \textit{minimal} (resp., \textit{maximal})
if $\DES(w)=\emptyset$ (resp., $\ASC(w)=\emptyset$), that is 
for each $i=1,2,\ldots,{n-1}$ the following property holds:\\[-0.7cm]
 \begin{align}
  \text{if $x_{i}$ and $x_{i+1}$ are distinct and comparable, then $x_{i}<x_{i+1}$ (resp., $x_{i}> x_{i+1}$)}.
 \end{align}
The following result is due to Foata and Randrianarivony~\cite{FoRa}.
\begin{prop}[Proposition 2.2, \cite{FoRa}]\label{prop:min-max word}
 Each  equivalence class in $L(X,\leq)$ contains one and only one  minimal (resp., maximal) word.
\end{prop}
Foata and Randrianarivony gave an explicit 
construction (see Section~2 in~\cite{FoRa}) of the (unique) bijection $\theta$ that sends each minimal word in 
$X^*$ onto the maximal word that belongs to the same  equivalence class. 
We now recall Theorem~\ref{thm:EquiAscDesA} which can be seen as an 
extension of Proposition~\ref{prop:min-max word} and is the key (and main) result of the paper.
\begin{thm}\label{thm:EquiAscDes}
For any poset $(X,\leq)$,  the  set-valued statistics  $\DES$  and  $\ASC$, 
are equidistributed over each equivalence class $\big[w\big]\in {L(X,\leq)}$, i.e.,
for any set $S$, there are as many words $w'$ in $\big[w\big]$ satisfying $\DES(w')=S$ as those satisfying $\ASC(w')=S$.
\end{thm}


\subsection{The partially commutative monoid $L(\B,\ll)$}

 Consider the poset $(\B,\ll)$ where $\B\subseteq \Pm$ is the set of (finite) decreasing sequences 
(or equivalenty, nonempty finite subsets) of positive integers and $\ll$ is the partial order on~$\Pm$
introduced in Definition~\ref{def:partial order}, that it is for $D,D'\in \B$, we have
$$
\text{$D \ll D' \iff $ $D=D'$ or $D$ is below $D'$}.
$$
 By convenience, the letters of a word in $\Bm$ 
will be separated by vertical lines, and the equivalence class in $L(\B,\ll)$ 
of an element $\pi\in\Bm$ will be denoted by $\big[\pi\big]$.  For example,
the word $\pi$ the letters of which are from left to right $6\,5\,3$, $2\,1$ and $3$ is written as 
$\pi=6\,5\,3\,|\,2\,1\,|\,3$ and its equivalence class $\big[\pi\big]\in L(\B,\ll)$ is 
\begin{align*}
 \big[\pi\big]=\{6\,5\,3\,|\,2\,1\,|\,3    \,,\,   2\,1\,|\,6\,5\,3\,|\,3    \,,\,   6\,5\,3\,|\,3\,|\,2\,1\}.
\end{align*}
Note that $\OP_n^k$ is just a subset of $\Bm$. In particular, we can extend the notion of \emph{block descent} 
(resp., \textit{block rise}) to $\Bm$. For $\pi=D_1\,|\,D_2\,|\cdots|\,D_k\in \Bm$, we set 
\begin{align*}
 \bDES(\pi)&=\{i\,:\,1\leq i\leq k-1\quad\text{and}\quad D_i\gg D_{i+1}\}\\
 \bASC(\pi)&=\{i\,:\,1\leq i\leq k-1\quad\text{and}\quad D_i\ll D_{i+1}\}.
\end{align*} 
As an immediate consequence of Theorem~\ref{thm:EquiAscDes}, we obtain the following result.
\begin{cor}\label{cor:EquiAscDes-Bm}
For any $\pi\in\Bm$, the 
set valued statistics $\bDES$ and $\bASC$ have the same
distribution over the equivalence class $\big[\pi\big]\in L(\B,\ll)$.
\end{cor}
 It will be convenient to extend the definition of pattern containment in $\Pm$ 
to $\Bm$ (and thus to $\OP_n^k$).
\begin{defn}\label{def:pattern in Bn}
Let $\pi=D_1\,|\,D_2\,|\,\cdots\,|\,D_k$ be an element of $\Bm$  and 
$p=\tau^{(1)}-\tau^{(2)}-\cdots-\tau^{(\ell)}$  be a piecewise decreasing pattern
of type $(j_1, j_2,\cdots, j_{\ell})$.

We say that a word $w^{(1)}\,|\,w^{(2)}\,|\,\cdots\,|\,w^{(\ell)}\in \Bm$ 
 is an occurrence of $p$ in~$\pi$ if there exist indices $1\leq t_1<t_2<\cdots<t_{\ell}\leq k$
such that 
\begin{itemize}
\item for $i=1,2,\ldots,\ell$, $w^{(i)}$ is a contiguous subsequence of $D_{t_i}$ of length $j_i$;
 \item the word $w^{(1)}\cdot w^{(2)}\cdot\cdots\cdot w^{(\ell)}$  is order isomorphic to 
$\tau^{(1)}\cdot \tau^{(2)}\cdot\cdots\cdot \tau^{(\ell)}$.
\end{itemize}
\end{defn}
For example, the word $5\,3\,2\,|\,6\,4\,1\,|\,5\,4$ does not contain the 
pattern $3\,1-4\,2-4$, but it contains exactly one occurrence of the pattern $3\,1-4\,2-3$, namely the
subword the letters of which are boldfaced  in ${\bf 5\,3}\,2\,|\,{\bf 6\,4}\,1\,|\,{\bf 5}\,4$.
If $p_1, p_2, \ldots, p_m$ are dashed patterns, we let $(p_1,\ldots, p_m)$ denote the $m$-statistic 
which associates to a word $\pi\in\Bm$ 
the vector $(a_1,a_2,\cdots,a_m)$ where $a_i$ is the number of occurences of  $p_i$ in $\pi$.
\begin{prop}\label{prop:pattern in Bn}
For any word $\pi$ in $\Bm$ and any piecewise decreasing and connected pattern $p$, 
the parameter $(p)$ is constant on the equivalence class $\big[\pi\big]\in L(\B,\ll)$.
\end{prop}

\pf Let $\pi$ and $\pi'$ be two adjacent words in $\Bm$; then, there exist $D_1,D_2,\ldots,D_k$ in $\B$ 
and a positive integer $j<k$ such that $D_j$ and $D_{j+1}$ are distinct and comparable 
(i.e., $D_j\gg D_{j+1}$ or $D_j\ll D_{j+1}$) and 
$$\pi=D_1\,|\,D_2\,|\,\cdots\,|\,D_k\quad\text{and}\quad
\pi'=D_1\,|\,D_2\,|\,\cdots\,|\,D_{j-1}\,|\,D_{j+1}\,|\,D_{j}\,|\,D_{j+2}\,|\cdots |\,D_k.$$

 Let $p=\tau^{(1)}-\tau^{(2)}-\cdots-\tau^{(\ell)}$  be a piecewise decreasing and connected dashed pattern.
We have to show that  $(p)(\pi)=(p)(\pi')$. Clearly, it suffices, by symmetry, to prove that  $(p)(\pi)\leq (p)(\pi')$.
Suppose that  $w^{(1)}\,|\,w^{(2)}\,|\,\cdots\,|\,w^{(\ell)}$ 
 is an occurrence of $p$ in~$\pi$ such that  $w^{(i)}$, $i=1,2,\ldots,\ell$, is a contiguous subsequence of 
$D_{t_i}$ for ${1\leq t_1<t_2<\cdots<t_{\ell}\leq k}$. 
\begin{enumerate}
\item If $|\{j,j+1\}\cap\{t_1,t_2,\ldots,t_{\ell}\}|\leq 1$, it is obvious to see that $w^{(1)}\,|\,w^{(2)}\,|\,\cdots\,|\,w^{(\ell)}$ 
is still an occurrence of $p$ in $\pi'$; thus $(p)(\pi)\leq (p)(\pi')$
\item We can not have $|\{j,j+1\}\cap\{t_1,t_2,\ldots,t_{\ell}\}|=2$. Suppose the contrary. Then, $j=t_m$ and $j+1=t_{m+1}$
for a certain integer $m$. Since the pattern $p$ is connected, the sequences $w^{(m)}$ and $w^{(m+1)}$ are equal or 
incomparable. But, $w^{(m)}$ and $w^{(m+1)}$ are contiguous subsequences of $D_j$ and $D_{j+1}$ respectively, whence 
$D_j$ and $D_{j+1}$ are equal or incomparable. This contradicts our assumption.
\end{enumerate}
Altogether, this implies that $(p)(\pi)\leq (p)(\pi')$.

\qed

Combining Corollary~\ref{cor:EquiAscDes-Bm} and Proposition~\ref{prop:pattern in Bn}, we arrive at the following 
result.
\begin{thm}\label{thm:key result-local}
For any word $\pi\in\Bm$ and any piecewise decreasing and connected patterns $p_1, p_2,\ldots,p_m$, the parameters 
$\big(\bDES,(p_1,p_2,\ldots,p_m)\big)$ and $\big(\bASC,(p_1,p_2,\ldots,p_m)\big)$ 
have the same distribution over the class $\big[\pi\big]\in L(\B,\ll)$.
\end{thm}

 Given a multiset $\M=\{D_1,D_2,\ldots,D_r\}\subset \B$, 
let $\R\big(\M\big)$ denote  the set of words in $\Bm$ of length~$r$ the letters of which 
are $D_1$, $D_2$,\ldots,$D_r$, that is
$$\R\big(\M\big)=\{ \,D_{\sigma(1)}\,|\, D_{\sigma(2)}\,|\,\cdots\,|\,D_{\sigma(r)}\,;\,\sigma\in\S_r\}.$$ 
For example, we have 
$\R\big(\{2\,1\,,\,2\,1\,,\,5\,3\}\big)=\{2\,1 \,|\, 2\,1  \,|\,  5\,3\,,\, 2\,1\,|\,5\,3\,|\,2\,1\,
,\,  5\,3\,|\,2\,1\,|\,2\,1\}.$
It is clear that for all $\pi\in\R(\M)$, we have $\big[\pi\big]\in\R\big(\M\big)$. The following result is
immediate from Theorem~\ref{thm:key result-local}.

\begin{cor}\label{cor:key result-local}
For any multiset $\M\subset \B$ and any piecewise decreasing and connected patterns $p_1, p_2,\ldots,p_m$, the parameters 
$\big(\bDES,(p_1,p_2,\ldots,p_m)\big)$ and $\big(\bASC,(p_1,p_2,\ldots,p_m)\big)$ 
have the same distribution over $\R\big(\M\big)$.
\end{cor}

\subsection{Applications}
\subsubsection{Wilf classification of dashed patterns: Proof of Theorem~\ref{thm:Wilf-equiv-local}}

It is often convenient  to identify a word $w=w_1w_2\ldots w_n\in\P^*$
with the element $\pi^w$ of $\Bm$ the letters of which are from left to right the descending runs 
of $w$, i.e,  $\pi^w$ is obtained from $w$ by placing a vertical line between 
$w_j$ and $w_{j+1}$ whenever $w_j\leq w_{j+1}$.  For example, we have the correspondence
$$
w=3\,5\,4\,1\,6\,5\,5\,3\,6\,5 \;  \longleftrightarrow \;\pi^w=3  \,|\, 5\,4\,1  \,|\, 6\,5  \,|\, 5\,3 \,|\, 6\,5. 
$$ 
 It is easy to see that the map $w\to\pi^w$ establishes a bijection between $\P^*$ and 
the set of minimal elements in $\Bm$ (i.e., those which have no block descent). 
In particular, by Proposition~\ref{prop:min-max word}, the map $w\to [\pi^w]$ establishes a bijection 
between $\P^*$ and $L(\B,\ll)$.

Let $p=\tau^{(1)}-\tau^{(2)}-\cdots-\tau^{(\ell)}$ be a piecewise decreasing pattern. 
By definition, in an occurrence of $p$ in a word $w\in\Pm$, the letters corresponding to a  ``bloc'' $\tau^{(i)}$ 
have to form a contiguous subsequence of a descending run of $w$. This leads to
\begin{fact}\label{fact:pattern in B vs W}
For any piecewise decreasing pattern $p$, the number of occurrences 
of $p$ in a word $w\in\Pm$ is equal to the number of occurrences  of $p$ in $\pi^w(\in\Bm)$ (in the sense 
of Definition~\ref{def:pattern in Bn}).
\end{fact}

In the remainder of this section, we suppose we are given $m$ decreasing and connected patterns $p_1, p_2,\ldots,p_m$.
 It is easy to check that, for all multisets $\M\subseteq \B$, 
the reverse map $\Rev$ that sends each 
word~$\pi=D_1\,|\,D_2|\cdots|\,D_k$ in $\R(\M)$ onto the 
word~$\Rev \pi=D_k\,|\,D_{k-1}\,|\cdots|\,D_1$ is an involution on $\R(\M)$
such that, for all $\pi\in \R(\M)$, we have
\begin{align}\label{eq:propRev}
 \big|\bDES(\Rev \pi)\big|= \big|\bASC(\pi)\big| \;,\;
(p_1,p_2,\ldots,p_m)(\Rev \pi)=(\Rev p_1,\Rev p_2,\ldots,\Rev p_m)(\pi).
\end{align}
Let $\R_{\emptyset}(\M)$ be the subset of $\R(\M)$ which consists of the minimal elements in $\R(\M)$, i.e.,
$$
\R_{\emptyset}(\M)=\{\pi\in\R(\M)\,:\,\bDES(\pi)=\emptyset\}.
$$
Recall that for a multiset $\M=\{D_1,D_2,\ldots,D_r\}\subset\B$, we let  
$\W(\M)$ denote the set of words in $\P^*$ the descending runs of which are $D_1$, $D_2$,\ldots,$D_r$.
The map $w\in\Pm\to\pi^w\in\Bm$ sends bijectively
$\W(\M)$ onto the set $\R_{\emptyset}(\M)$.
For example, if $\M=\{4\,2\,1   \,,\,  6\,5  \,,\,  7\,5\}$, the reader can check that  
$$\W(\M)=\{ 4\,2\,1\,6\,5\,7\,5\,  \,,\,  4\,2\,1\,7\,5\,6\,5 \}\quad\text{and}\quad
\R_\emptyset(\M)=\{4\,2\,1   \,|\,  6\,5 \,|\, 7\,5  \,,\,  4\,2\,1 \,|\, \,7\,5 \,|\, \,6\,5 \,\}
$$
Use of Fact~\ref{fact:pattern in B vs W}, \eqref{eq:propRev} and Corollary~\ref{cor:key result-local},
we see that the parameters $(p_1,p_2,\ldots,p_m)$ and ${(\Rev p_1,\Rev p_2,\ldots,\Rev p_m)}$
are equidistributed over $\R_{\emptyset}(\M)$, and thus over, $\W(\M)$. This ends the proof 
of Theorem~\ref{thm:Wilf-equiv-local} which refines Theorem~\ref{thm:Wilf-equiv}.

\subsubsection{Statistics on ordered set partitions}

\paragraph{Proof of Theorem~\ref{thm:application1}.} 

It is easy to check that for any ordered set partition  $\pi$, 
the value $\rsb \pi$ (resp., $\lsb \pi$) is equal to the number of occurrences, 
in the sense of Definition~\ref{def:pattern in Bn}, of the pattern $2-3\,1$ (resp., $3\,1-2$) in $\pi$. 
\begin{fact}\label{fact:OPstat-pattern}
For any ordered set partition $\pi\in\OP_n^k$, we have 
$\rsb \pi=(2-3\,1)(\pi)$ and $\lsb \pi=(3\,1-2)(\pi)$.
\end{fact}
 Since $\OP_n^k\subset\Bm$ and any rearrangement of the blocks
of a partition $\pi\in\OP_n^k$ is still in $\OP_n^k$, we see that the equivalence 
class $\big[\pi\big]\in L(\B,\ll)$ is contained 
in $\OP_n^k$. 
For example, if $\pi=5\,3\,|\,4\,|\,2\,1\in \OP_5^3$, we have
\begin{align*}
 \big[\pi\big]=\{5\,3\,|\,4\,|\,2\,1 \,,\, 5\,3\,|\,2\,1\,|\,4 \,,\, 2\,1\,|\,5\,3\,|\,4\}\subseteq \OP_5^3.
\end{align*}
 Moreover, it is obvious that the set valued statistics
$\open$ and $\clos$ are constant on each equivalence class~$\big[\pi\big]$.
This, combined with Fact~\ref{fact:OPstat-pattern} and Theorem~\ref{thm:key result-local},
leads immediately to the following result.
\begin{cor}\label{cor:key-OP}
For any ordered set partition $\pi\in\OP_n^k$, the 
parameters $$(\bDES,\open,\clos,\rsb,\lsb)\quad\text{and}\quad(\bASC,\open,\clos,\rsb,\lsb)$$ 
have the same 
distribution over the equivalence class~$\big[\pi\big]$, and thus over $\OP_n^k$.
\end{cor}
 
It follows immediately from the definition of the parameters $\MAK$ and ${\MAK}'$ and  
Corollary~\ref{cor:key-OP} that the 
3-statistics ${(\MAK,{\MAK}',\bDES)}$ and $(\MAK,{\MAK}',\bASC)$ have the same
distribution  over $\OP_n^k$.
   Consider the transformation \textit{complement}~$\c$ that maps each ordered set 
partition~$\pi$ of~$\Nn$ onto the ordered set partition~$\c\pi$ of~$\Nn$
obtained from~$\pi$ by complementing each of the letters in~$\pi$,
 that is, by replacing $i$ by $n+1-i$. For example, if $\pi=5\,3\,|\,4\,|\,2\,1$, 
then $\c \pi=1\,3\,|\,2\,|\,4\,5\equiv 3\,1\,|\,2\,|\,5\,4$.  It is easily checked that the
 complement~$\c$ sends the parameter~$(\MAK,{\MAK}',\bASC)$ onto~$({\MAK}',\MAK,\bDES)$. 
Since the complement map is an involution on $\OP_n^k$,
we obtain that  $(\MAK,{\MAK}',\bDES)$ and $({\MAK}',\MAK,\bDES)$ are equidistributed over $\OP_n^k$.
This concludes the proof of Theorem~\ref{thm:application1}.

\paragraph{A new Euler-Mahonian statistic.}
 It is worth noting that Corollary~\ref{cor:key-OP} permits to get 
new hard Euler-Mahonian statistics from known ones. For example,
using the following result of Zeng and the author~\cite{KaZe}\\[0.2cm]
\textbf{Result C}. \textit{For $n\geq k\geq 1$, the partition statistic $\lsb-\bMAJ+\,k(k-1)$ 
is Euler-Mahonian on~$\OP_n^k$};\\[0.2cm]
and elementary properties of the reverse map $\Rev$,
we obtain the following result.
\begin{thm}\label{thm:application2}
Given a partition $\pi=B_1\,|\,B_2\,|\cdots|\,B_k$, 
let $\nbDes(\pi)$ be the number of block non-descents of $\pi$, i.e. $\nbDes(\pi)=k-1-|\bDES(\pi)|$.
Then, the partition statistic $\STAT:=\rsb+k\cdot \nbDes+\bMAJ$ is Euler-Mahonian on $\OP_n^k$,
\end{thm}
\pf By result~C, it suffices to prove that the parameter $\STAT$ is 
equidistributed with the parameter $\lsb-\bMAJ+\,k(k-1)$ over $\OP_n^k$. It is easy to check that 
reverse map  $\Rev$ is an involutive transformation of $\OP_n^k$ and 
sends the parameter ${(\lsb,\bASC)}$ onto the parameter~$(\rsb,k-\bDES)$.
This, combined with Corollary~\ref{cor:key-OP},  implies that 
the parameters~$(\lsb,\bDES)$ and~$(\rsb,k-\bDES)$, and thus the parameters $\STAT=\rsb-(k\cdot \bDes-\bMAJ)+k(k-1)$ 
and $\STAT'=\lsb-\bMAJ+k(k-1)$, have the same distribution over $\OP_n^k$. 
\qed

%
%

\section{Descent and rise sets in equivalence classes of a partially commutative monoid}

  This section is dedicated to the proof of Theorem~\ref{thm:EquiAscDesA} which asserts that the  
set-valued statistics $\DES$ and~$\ASC$ are equidistributed over each equivalence class $\big[w\big]\in L(X,\leq)$ 
for any nonempty poset $(X,\leq)$. We shall give an ``algebraic proof'' as well as a bijective one.

\subsection{A proof by the inclusion-exclusion principle}

Let $w$ be a word in $X^*$. For a set $S\subseteq \P$, we set 
\begin{align*}
 f_{=}(S):=|\{w' \in \big[w\big] \,: \;\DES(w')=S\}|\quad \text{and}\quad
 g_{=}(S):=|\{w' \in \big[w\big]\, : \;\ASC(w')=S\}|.
\end{align*}
We have to show that $f_{=}(S)=g_{=}(S)$. For a set $T\subseteq \P$, define
\begin{align*}
f_{\subseteq}(T):=|\{w' \in \big[w\big]\, :\; \DES(w')\subseteq T\}|\quad \text{and}\quad
g_{\subseteq}(T):=|\{w' \in \big[w\big]\, :\; \ASC(w') \subseteq T\}|,
\end{align*}
so that we have $f_{\subseteq}(S)=\sum_{T\subseteq S} f_{=}(T)$ and $g_{\subseteq}(S)=\sum_{T\subseteq S} g_{=}(T)$. 
Then, by the principle of inclusion-exclusion (see e.g. Chapter~5 in~\cite{Aigner}), we have
\begin{align}\label{eq:incl-excl}
 f{=}(S)=\sum_{T\subseteq S}(-1)^{|S|-|T|} f_{\subseteq}(T) \quad\text{and}\quad
 g{=}(S)=\sum_{T\subseteq S}(-1)^{|S|-|T|} g_{\subseteq}(T).
\end{align}
Suppose that  $T=\{1\leq t_1<t_2<\cdots<t_{k-1}<n\}$. A moment's thought will convince the reader 
that $f_{\subseteq}(T)$ (resp., $g_{\subseteq}(T)$) is exactly the number of $k$-tuples
$(w^{(1)},w^{(2)},\cdots,w^{(k)})$ in $\big(\Pm\big)^k$ such that 
\begin{itemize}
 \item the word $w'=w^{(1)}\cdot w^{(2)}\cdots w^{(k)}$ obtained by concatenation of 
the $w^{(i)}$'s belongs to the class $\big[w\big]$;
 \item $|w^{(i)}|=t_i-t_{i-1}$ for $1\leq i\leq k$, with $t_0=0$ and $t_{k}=n$,  where $n$ is the length of $w$  and 
$|w^{(i)}|$ is the length of $w^{(i)}$;
 \item $w^{(i)}$ is minimal (resp., maximal) for $1\leq i\leq k$.
\end{itemize}
Recall that $\theta$ denotes the (unique) bijection that sends each minimal word in $X^*$ onto the maximal word
that belongs to the same  equivalence class. Then, it is easy to see that the  bijection
that sends each $k$-tuple of  minimal words $(w^{(1)},w^{(2)},\cdots,w^{(k)})$ onto the $k$-tuple of  maximal
words~$( \theta\,w^{(1)},\theta\,w^{(2)},\cdots,\theta\,w^{(k)})$ leads
to the identity~$f_{\subseteq}(T)=g_{\subseteq}(T)$, from which we deduce, by~\eqref{eq:incl-excl}, that 
$f_{=}(S)=g_{=}(S)$. This concludes the proof of Theorem~\ref{thm:EquiAscDes}.

\subsection{A bijective proof}

\subsubsection{An extremal case: $(X,\leq)$ is a totally ordered set} 

Suppose $(X,\leq)$ is a totally ordered set. Without loss of generality, we can assume
$X=\{1,2,\ldots,r\}$, the set of the $r$ first  positive integers, with the natural order
$1<2<\cdots<r$. In this case, the equivalence class of a word $w\in X^*$ 
is just its rearrangement class (i.e., the set of words obtained by permuting the letters
of $w$). Denote by $\R(n_1,n_2,\ldots,n_r)$ the rearrangement class of the word $1^{n_1}2^{n_2}\ldots r^{n_r}$,
or equivalently  the set of words which have exactly $n_1$ occurrences of the ''letter'' 1, 
$n_2$ occurences of 2,\ldots, $n_r$ occurences of $r$. In this particular case, Theorem~\ref{thm:EquiAscDesA} just asserts
that for any $r$-tuple $(n_1,n_2,\ldots,n_r)$ of nonnegative integers,
the  set-valued statistics $\DES$ and $\ASC$ are equidistributed over $\R(n_1,n_2,\ldots,n_r)$,
or by abuse of notation,
\begin{align}\label{eq:EquiAscDes_multipermutations}
 \sum_{w\in \R(n_1,n_2,\ldots,n_r)}q^{\DES(w)}= \sum_{w\in \R(n_1,n_2,\ldots,n_r)}q^{\ASC(w)}.
\end{align}
It is not hard to give a bijective proof of~\eqref{eq:EquiAscDes_multipermutations}. Indeed, 
for $i=1,2,\ldots,r-1$, there is a simple bijection 
$\gamma_i:\R(n_1,n_2,\ldots,n_r)\mapsto \R(n_1,n_2,\ldots n_{i-1},n_{i+1},n_{i},n_{i+2},\ldots,n_r)$
which preserves the descent set (see e.g. the proof of Theorem~10.2.1 in~\cite{Lo}):
``consider a word $w\in\R(n_1,n_2,\ldots,n_r)$ and write all its factors of the form $(i+1)i$ 
in bold-face; then replace all the maximal factors $i^a(i+1)^b$ with $a\geq0$, $b\geq0$, that do not involve 
any bold-face letters by $i^b(i+1)^a$. Finally, rewrite all the bold-face letters in roman type.'' It follows that 
the map ${\mathbf \rho}=\gamma_{1}\;\gamma_2\gamma_1\;\cdots\;\gamma_{r-2}\cdots\gamma_2\gamma_1\;
\gamma_{r-1}\cdots\gamma_2\gamma_1$ establishes a bijection from $\R(n_1,n_2,\ldots,n_r)$ to 
$\R(n_r,n_{r-1},\ldots,n_1)$ which preserves the descent set. Let $\c$ be the involutive transformation complement
that maps each word $w=w_1w_2\ldots w_n$ 
onto $\textbf{c}\,w=(r+1-w_1)(r+1-w_2)\ldots(r+1-w_n)$. Then, it is easily checked that the map
$\c\circ\rho:\R(n_1,n_2,\ldots,n_r)\mapsto \R(n_1,n_2,\ldots,n_r)$ is a bijection which sends
the descent set onto the rise set. 

\subsubsection{The general case: a bijection based on the involution principle}
 The general involution principle, 
introduced by Garsia and Milne~\cite{GaMi}, permits to convert non-bijective proofs (notably, 
proofs by inclusion-exclusion principle) into bijective ones. 

Let $(X,\leq)$ be a nonempty poset. For any (finite) subset $S\subseteq \P$, let $D_S$ and $A_S$ be the subsets of $X^*$ defined by
\begin{align*}
 D_S:=\{w \in X^*\,: \;\DES(w)=S\}\quad\text{and}\quad A_S:=\{w \in X^*\,: \;\ASC(w)=S\}.
\end{align*}
Using the general involution principle, we construct a bijection $\Gamma_S:D_S\mapsto A_S$ for any (finite) 
subset $S\subseteq \P$. In order to invoke the general involution principle 
 (see e.g.~\cite{GaMi} or Chapter~5 in~\cite{Aigner}), we need 
(we assume that the reader is familiar with the used terminology) 
\begin{itemize}
 \item two signed sets $Y=Y^{+}\sqcup Y^{-}$ and $Z=Z^{+}\sqcup Z^{-}$, with $D_S\subseteq Y^{+}$
and $A_S\subseteq Z^{+}$;\\[-0.7cm]
 \item an alternating involution $\psi$ on $Y=Y^{+}\sqcup Y^{-}$ with $\Fix(\psi)=D_S$, where 
$\Fix(\psi)$ stands for the set of fixed points of $\psi$;\\[-0.7cm]
 \item an alternating involution $\phi$ on $Z=Z^{+}\sqcup Z^{-}$ with $\Fix(\phi)=A_S$;\\[-0.7cm]
 \item a sign-preserving bijection $F$ from $Y$ to $Z$. \\[-0.7cm]
\end{itemize}

\paragraph{The signed sets.}
Let $Y=Y^{+}\sqcup Y^{-}$ be the signed set 
$$Y=\{(w,T) \,: \;w \in X^*,\;S\subseteq T\subseteq \DES(w)\}$$
with $(w,T)\in Y^{+}$ if $|T|-|S|$ is even, and $(w,T)\in Y^{-}$ if $|T|-|S|$ is 
odd. By identifying an element $w\in D_S$ with the pair $(w,S)\in Y^{+}$, we have 
$D_S\subseteq Y^{+}$.

Similarly, let $Z=Z^{+}\sqcup Z^{-}$ be the signed set 
$$Z=\{(w,T) \,: \;w \in X^*,\; S\subseteq T\subseteq \ASC(w)\}$$
with $(w,T)\in Z^{+}$ if $|T|-|S|$ is even, and $(w,T)\in Z^{-}$ if $|T|-|S|$ is 
odd. By identifying an element $w\in A_S$ with the pair $(w,S)\in Z^{+}$, we have 
$A_S\subseteq Z^{+}$.\\[-0.7cm]

\paragraph{The alternating involution $\phi$ on $Y=Y^{+}\sqcup Y^{-}$.} 
For $(w,T)\in Y$, let $d(w)=\max(\DES(w)\setminus S)$ whenever $\DES(w)\neq S$.
Then define $\phi:Y\mapsto Y$ by
$$
\phi(w,T)=\left\{
    \begin{array}{ll}
     (w,T\setminus\{d(w)\}), & \hbox{if $d(w)\in T$,} \\
     (w,T\cup\{d(w)\}), & \hbox{if $d(w)\notin T$,}
     \end{array}
\right. 
$$
whenever $\DES(w)\neq S$, and $\phi(w,S)=(w,S)$ if $\DES(w)= S$.
The mapping $\phi$ is clearly an alternating involution on $Y=Y^{+}\sqcup Y^{-}$, and 
the only fixed points are the elements $(w,S)$ with $\DES(w)=S$ which were 
identified with the elements of $D_S$.\\[-0.7cm]

\paragraph{The alternating involution $\psi$ on $Z=Z^{+}\sqcup Z^{-}$.}
For $(w,T)\in Z$, let $a(w)=\max(\ASC(w)\setminus S)$ whenever $\ASC(w)\neq S$.
Then define $\psi:Z\mapsto Z$ by
$$
\psi(w,T)=\left\{
    \begin{array}{ll}
     (w,T\setminus\{a(w')\}), & \hbox{if $a(w)\in T$,}\\
     (w,T\cup\{a(w')\}), & \hbox{if $a(w)\notin T$,}
    \end{array}
\right. 
$$
whenever $\ASC(w)\neq S$, and $\psi(w,S)=(w,S)$ for $\ASC(w)= S$.
The mapping $\psi$ is clearly an alternating involution on $Z$, and the only fixed points
are the elements $(w,S)$ with $\ASC(w)=S$ which were identified with the elements of $A_S$.\\[-0.7cm]

\paragraph{The sign-preserving bijection $F$ from $Y$ to $Z$.}
Recall that, for $T\subseteq \P$, the $T$-factorization of a word $w=x_1x_2\cdots x_n$ in $X^*$
is the (unique) factorization $w=w^{(1)}\cdot w^{(2)}\cdot \ldots\cdot w^{(k)}$ of $w$ defined by
$x_i$ is the last letter of a factor $w^{(j)}$ if and only if $i\notin T$ or $i=n$.
Then define $F:Y\mapsto Z$ as follows: if $(w,T)\in Y$ and the word $w$ has $T$-factorization 
$w=w^{(1)}\cdot w^{(2)}\cdots w^{(k)}$,
we set $F(w,T)=({\bf r}\,w^{(1)}\cdot {\bf r}\,w^{(2)}\cdots {\bf r}\,w^{(k)},T)$, 
where $\r$ is, as usual, the reverse image. 

In order to see that $F$ is well defined from $Y$ to $Z$, it suffices 
(since $S\subseteq T$) to prove that $T\subseteq \ASC({\bf r}\,w^{(1)}\cdot {\bf r}\,w^{(2)}\cdots {\bf r}\,w^{(k)})$.
By definition of $Y$, we have $T\subseteq \DES(w)$ and thus, by definition of the $T$-factorization,
the factors $w^{(j)}$, $1\leq j\leq k$, are decreasing words. This implies that 
the words $\r\,w^{(j)}$, $1\leq j\leq k$, are increasing words, from which it is immediate to see that 
$T\subseteq \ASC({\bf r}\,w^{(1)}\cdot {\bf r}\,w^{(2)}\cdots {\bf r}\,w^{(k)})$. Moreover, since the letters 
of an increasing (or decreasing) word are pairwise comparable, the words $w^{(j)}$ and $\r\,w^{(j)}$, $1\leq j\leq k$,
are equivalent, whence ${\bf r}\,w^{(1)}\cdot {\bf r}\,w^{(2)}\cdots {\bf r}\,w^{(k)}\in\big[w\big]$.
 It is clear that $F$ is sign-preserving. To see that $F$ is bijective, we construct its inverse.
Define $G:Z\mapsto Y$ as follows: if $(w,T)\in Z$ and the word $w$ has $T$-factorization 
$w=w^{(1)}\cdot w^{(2)}\cdots w^{(k)}$, we set $G(w,T)=({\bf r}\,w^{(1)}\cdot {\bf r}\,w^{(2)}\cdots {\bf r}\,w^{(k)},T)$. 
It is easy to check that $G$ is well-defined and that $G=F^{-1}$. Summarizing, we have obtained the following result.

\begin{prop}
 The mapping $F:Y\to Z$ is a sign-preserving bijection such that  
if $(w,T)\in Y$ and $F(w,T)=(w',T)$, then $w$ and $w'$ are equivalent, 
i.e., $w'\in \big[w\big]$. 
\end{prop}

Applying the general involution principle (see e.g. Chapter~5 in~\cite{Aigner}), 
we arrive at
\begin{thm}\label{thm:bijDESASC}
Let $w\equiv(w,S)\in D_S$, i.e., $w\in X^*$ and $\DES(w)=S$. There is a least integer $o(w)$ such that 
$$
F\left(\phi F^{-1} \psi F\right)^{o(w)}(w,S)\in A_S.
$$
If we set $(w',S):=F\left(\phi F^{-1} \psi F\right)^{o(w)}(w,S)$, then 
the mapping $\Gamma_S: w\mapsto w'$ is a bijection from $D_S$ to $A_S$
such that $w'\in\big[w\big]$. 

If we by $\Gamma$ the transformation of $X^*$
whose restriction on $D_S$ is equal to $\Gamma_S$, then 
$\Gamma$ is a bijective transformation of $X^*$ such that 
$\DES(w)=\ASC(\Gamma\,w)$ and $\Gamma\,w\in\big[w\big]$ for all $w\in X^*$. 
\end{thm}

 It follows from Proposition~\ref{prop:min-max word} that there exists a unique bijection 
$\theta:D_{\emptyset}\mapsto A_{\emptyset}$ which satisfies the following condition:
\begin{align*}
 \text{for all $w\in D_{\emptyset}$, $\theta(w)$ and $w$ are equivalent, i.e. $\theta(w)\in \big[w\big]$.}
\end{align*}
Therefore, by Theorem~\ref{thm:bijDESASC}, $\theta$ is the restriction of $\Gamma$ on $\D_{\emptyset}$, 
i.e., $\theta=\Gamma_{\emptyset}$. 
Foata and Randriarivony~\cite{FoRa} gave a more direct description of the mapping $\theta$, which 
can be defined recursively as follows (see~\cite[Section 6]{ClStZe}).
\paragraph{The mapping $\theta:D_{\emptyset}\mapsto A_{\emptyset}$.} 
\begin{itemize}
 \item if $|w|\leq 1$, $\theta(w)=w$; 
 \item if $|w|=n\geq 2$: suppose $w=w'\cdot x$ with $w'=x_1 x_2\cdots x_{n-1}$ and $x$ in $X$. Let $t$ be 
the largest integer $\leq n-1$ such that $x_t$ is incomparable and distinct from $x$. Then $\theta(w)$ is obtained from $\theta(w')$ 
by inserting $x$ between $x_t$ and $x_{t+1}$. If there is no such~$t$, then we set $\theta(w)=x\cdot \theta(w')$.
\end{itemize}
Illustrations are given in the next subsection. 


\subsubsection{Illustration: the transformation $\Gamma:\Bm\mapsto \Bm$ and Clarke et al.'s 
                                                    transformation ${\epsilon:\Pm\to \Pm}$}

 Theorem~\ref{thm:bijDESASC} leads to the following result. 
\begin{thm}\label{thm:bij-gamma}
Consider the monoid $\Bm$ generated by $(\B,\ll)$. Then, the transformation 
 $\Gamma:\Bm\mapsto \Bm$ described in Theorem~\ref{thm:bijDESASC} is a bijection which sends the parameter $\bDES$ onto $\bASC$,
and such that, for all $w\in \Bm$, $\Gamma\,w$ belongs to the equivalence class $\big[w\big]\in L(\B,\ll)$.
\end{thm}

For example, if $w=3\,|\,9\,6\,|\,5\,4\,|\,2\,1\,|\,8\,7$, then $S=\bDES(w)=\{2,3\}$. 
Running the algorithm for $\Gamma_S$, we obtain 
$$ w\equiv(w,S)=(3\,|\,9\,6\,|\,5\,4\,|\,2\,1\,|\,8\,7,\{2,3\})
\overset{F}{\longmapsto} (3\,|\,2\,1\,|\,5\,4\,|\,9\,6\,|\,8\,7,\{2,3\})\in A_S,$$
whence $\Gamma(3\,|\,9\,6\,|\,5\,4\,|\,2\,1\,|\,8\,7)=3\,|\,2\,1\,|\,5\,4\,|\,9\,6\,|\,8\,7$.

If $w=2\,1\,|\,9\,6\,|\,5\,4\,|\,3\,|\,8\,7$, then $S=\bDES(w)=\{2,3\}$. 
Running the algorithm for $\Gamma_S$, we obtain 
\begin{align*}
&  w\equiv(w,S)=(2\,1\,|\,9\,6\,|\,5\,4\,|\,3\,|\,8\,7,\{2,3\})
\overset{F}{\longmapsto} (2\,1\,|\,3\,|\,5\,4\,|\,9\,6\,|\,8\,7,\{2,3\})\notin A_S\\
&\overset{\psi}{\longmapsto} (2\,1\,|\,3\,|\,5\,4\,|\,9\,6\,|\,8\,7,\{1,2,3\})
\overset{F^{-1}}{\longmapsto}(9\,6\,|\,5\,4\,|\,3\,|\,2\,1\,|\,8\,7,\{1,2,3\})\\
&\overset{\phi}{\longmapsto} (9\,6\,|\,5\,4\,|\,3\,|\,2\,1\,|\,8\,7,\{2,3\})
\overset{F}{\longmapsto}(9\,6\,|\,2\,1\,|\,3\,|\,5\,4\,|\,8\,7,\{2,3\})\notin A_S\\
&\overset{\psi}{\longmapsto} (9\,6\,|\,2\,1\,|\,3\,|\,5\,4\,|\,8\,7,\{2,3,4\})
\overset{F^{-1}}{\longmapsto} (9\,6\,|\,8\,7\,|\,5\,4\,|\,3\,|\,2\,1,\{2,3,4\})\\
&\overset{\phi}{\longmapsto} (9\,6\,|\,8\,7\,|\,5\,4\,|\,3\,|\,2\,1,\{2,3\})
\overset{F}{\longmapsto} (9\,6\,|\,3\,|\,5\,4\,|\,8\,7\,|\,2\,1,\{2,3\})\in A_S,
\end{align*}
whence $\Gamma(2\,1\,|\,9\,6\,|\,5\,4\,|\,3\,|\,8\,7)=9\,6\,|\,3\,|\,5\,4\,|\,8\,7\,|\,2\,1$.

If $w=3\,1\,|\,5\,4\,2\,|\,7\,6$, then $S=\bDES(w)=\emptyset$. 
Running the algorithm for $\Gamma_S$, we obtain 
\begin{align*}
&  w\equiv(w,S)=(3\,1\,|\,5\,4\,2\,|\,7\,6,\emptyset)
\overset{F}{\longmapsto} (3\,1\,|\,5\,4\,2\,|\,7\,6,\emptyset)\notin A_S
\overset{\psi}{\longmapsto} (3\,1\,|\,5\,4\,2\,|\,7\,6,\{2\})\\
&\overset{F^{-1}}{\longmapsto}(3\,1\,|\,7\,6\,|\,5\,4\,2,\{2\})
\overset{\phi}{\longmapsto} (3\,1\,|\,7\,6\,|\,5\,4\,2,\emptyset)
\overset{F}{\longmapsto}(3\,1\,|\,7\,6\,|\,5\,4\,2,\emptyset)\notin A_S\\
&\overset{\psi}{\longmapsto} (3\,1\,|\,7\,6\,|\,5\,4\,2,\{1\})
\overset{F^{-1}}{\longmapsto} (7\,6\,|\,3\,1\,|\,5\,4\,2,\{1\})
\overset{\phi}{\longmapsto} (7\,6\,|\,3\,1\,|\,5\,4\,2,\emptyset)\\
&\overset{F}{\longmapsto} (7\,6\,|\,3\,1\,|\,5\,4\,2,\emptyset)\in A_S,
\end{align*}
whence $\Gamma(3\,1\,|\,5\,4\,2\,|\,7\,6)=\theta(3\,1\,|\,5\,4\,2\,|\,7\,6)=7\,6\,|\,3\,1\,|\,5\,4\,2$. 
When $\bDes(w)=\emptyset$, it is more convenient to use the recursive description of $\theta$. For instance, 
if $w$ is as above, i.e. $w=3\,1\,|\,5\,4\,2\,|\,7\,6$, we have $\theta(3\,1\,|\,5\,4\,2)=3\,1\,|\,5\,4\,2$
and  
\begin{align*}
&  \Gamma(w)=\theta(w)=\theta(3\,1\,|\,5\,4\,2\,|\,7\,6)
=7\,6\,|\,3\,1\,|\,5\,4\,2.
\end{align*}

\paragraph{Clarke et al.'s involution ${\ep:\Pm\to \Pm}$.} 
 This transformation (see~\cite[Section 6]{ClStZe})
is a variation of the transformation $\theta=\Gamma_{\emptyset}$. 
If $w=x_1x_2\ldots x_n\in\P^*$, then $\ep(w)\in\P^*$ is the word obtained from 
$\Rev\big(\Gamma_{\emptyset}(\pi^w)\big)=\Rev\big(\theta(\pi^w)\big)$
by deleting the vertical bars. For example, if $w=3\,6\,4\,5\,3\,5\,3\,1\,7\,6$,
we have 
$$
\pi^w=3\,|\,6\,4\,|\,5\,3\,|\,5\,3\,1\,|\,7\,6
\overset{\Gamma_{\emptyset}}{\longmapsto} 6\,4\,|\,7\,6\,|\,3\,|\,5\,3\,|\,5\,3\,1 
\overset{\Rev}{\longmapsto} 5\,3\,1\,|\,5\,3\,|\,3\,|\,7\,6\,|\,6\,4
$$ 
whence $\ep(3\,6\,4\,5\,3\,5\,3\,1\,7\,6)=5\,3\,1\,5\,3\,3\,7\,6\,6\,4$.
Combining Theorem~\ref{thm:bij-gamma}, Proposition~\ref{prop:pattern in Bn} and Fact~\ref{fact:pattern in B vs W}
leads to the following result. 

\begin{thm}\label{thm:epsilon}
The map $\ep$ is a bijective transformation of $\Pm$ such that for any piecewise decreasing and connected patterns $p_1, p_2,\ldots,p_m$, 
we have $$(p_1,p_2,\ldots,p_m)(w)=(\Rev p_1,\Rev p_2,\ldots,\Rev p_m)\big(\epsilon(w)\big)$$
and $\big[\pi^{w}\big]=\big[\pi^{\ep (w)}\big]$ in $L(\B,\ll)$.
\end{thm}

It is worth noting that it is easy to give bijective proofs of Theorem~\ref{thm:application1} 
and Theorem~\ref{thm:application2} by using the map $\Gamma$.

%
%
\section{Concluding remarks and open problems}

\subsection{Wilf-classification of dashed patterns in subclasses of permutations}

In this paper, we focused our attention primarily on the (run- and reverse-complete) collections of permutations, 
of words and of compositions.
Other interesting run-complete (but not reverse-complete) collections are:\\[-0.7cm]
\begin{itemize}
\item the collection $\A^{(2)}=\bigcup_{n\geq0}\A_{2n}$ of reverse alternating permutations of even length, 
(and more generally, for any positive integer $k$, the collection $\A^{(k)}$ of permutations the descending runs of 
which are of length $k$);\\[-0.8cm]
\item the collection $\D^{(1)}=\bigcup_{n\geq0}\D^{(1)}_{n}$ (resp., $\D^{(3)}=\bigcup_{n\geq0}\D^{(3)}_{n}$) 
 of Dumont permutations of the first (resp., third) kind, where $\D^{(i)}_{n}$ is the set of Dumont permutations 
 of the $i$-th kind of length $n$ (see e.g.~\cite{BuJoSt} for a precise definition).\\[-0.8cm]
\end{itemize}
Note that there has been recent interest in the study of patterns in 
alternating permutations and Dumont permutations (see e.g.~\cite{Le,BuJoSt}).
Theorem~\ref{thm:Wilf-equiv} leads to wilf-classification results for these classes. For instance,
as an immediate consequence of Theorem~\ref{thm:Wilf-equiv}, we see that the joint
 distribution of the pair $((2-3\,1),(3\,1-2))$ over $\D^{(1)}_{n}$ is symmetric, which is not obvious since
$\D^{(1)}_{n}$ has no apparent symmetries (see the remark after Corollary~4.5 in~\cite{BuJoSt}). 
Note that Clarke et al.'s bijection leads to a simple bijective proof of this symmetry.

\subsection{Euler-Mahonian partition statistics and permutation statistics}

In this paper, we have shown that the equidistributions of the permutation statistics $(\des,\mak)$ and $(\des,\mak')$ 
over the symmetric group, and on the other part, of the statistics $\MAK+\bMAJ$ and ${\MAK}'+\bMAJ$ over ordered set partitions
have a very natural generalization in the context of ordered set partition statistics. 

It seems that there is a strong connection between the so called Euler-Mahonian 
permutations statistics and Euler-Mahonian partitions statistics, but the link is far to be understood. 
As an illustration, we present a conjecture which generalizes an important equidistribution result over the symmetric group
as well as an equidistribution result over set partitions.

Recall that a descent in a permutation $\sigma\in\S_n$ is 
an integer $i$, $1\leq i\leq n-1$, such that $\sigma(i)>\sigma(i+1)$. The \textit{Major index}, denoted $\maj$, 
is a well known Mahonian  permutation statistic. For a permutation $\sigma\in\S$, $\maj \sigma$ is defined as the 
sum of the descents in $\sigma$.
For example, if $\sigma=3\,2\,1\,7\,5\,6\,4\in\S_7$, then $1,2,4,6$ are the descents of $\sigma$ 
whence $\maj\sigma=1+2+4+6=13$. The following result is due to Foata and Zeilberger~\cite{FoZe}.\\[-0.3cm]

\textbf{Result C}. \textit{For $n\geq 1$, the bistatistics $(\des,\maj)$ and $(\des,\mak)$ 
are equidistributed over~$\S_n$}.\\[-0.3cm]

The Milne's statistic $\MIL$ is defined for an ordered set partition $\pi=B_1\,|\,B_2\,|\cdots|\,B_k\in\OP_n^k$
by $\MIL\,\pi=|B_2|+2|B_3|+\cdots+(k-1)|B_k|$. 
For example, $\MIL(8\,5\,|\,1\,|\,9\,6\,2\,|\,7\,4\,|\,3)=1\cdot1+2\cdot3+3\cdot2+4\cdot1=17$. 
Steingr\'{\i}msson~\cite{St} observed that the 
\textit{Major index} is closely related to the partition statistic~$\MIL$. More precisely,  he showed that if $\sigma\in\S_n$ has descending runs
$D_1,D_2,\ldots,D_k$ (i.e., $\pi^{\sigma}=D_1\,|\,D_2\,|\cdots|\,D_k$), then we have
\begin{align}\label{eq:MIL-maj}
\maj \sigma&=\MIL(D_1\,|\,D_2\,|\cdots|\,D_k)+{n+1\choose 2}-kn.
\end{align}
Steingr\'{\i}msson~\cite{St} proved that the 
statistic $\bmajMIL:=\MIL+\bMAJ$ is Euler-Mahonian on~$\OP_n^k$, while Zeng and the author~\cite{KaZe}
proved that the statistics $\MAK+\bMAJ$ (and ${\MAK}'+\bMAJ$) are Euler-Mahonian on~$\OP_n^k$.
Therefore, we have:\\[-0.3cm]

\textbf{Result D}. \textit{For $n\geq k\geq 1$, the partition statistics $\MIL+\bMAJ$ and $\MAK+\bMAJ$ 
are equidistributed over~$\OP_n^k$}.\\[-0.3cm]

  We suspect that actually more is true. The following conjecture has been verified for $n\leq 11$
with the help of Einar Steingr\'{\i}msson. 
\begin{conj}\label{conj:1}
 For $n\geq k\geq 1$, the bistatistics 
$$(\bDes,\MIL+\bMAJ)\quad\text{and}\quad (\bDes,\MAK+\bMAJ)$$
are equidistributed over $\OP_n^k$.
\end{conj}

This conjecture is interesting since it permits, in view of~\eqref{eq:MIL-maj} and~\eqref{eq:MAK-mak}, to unify Result~C and Result~D, and its solution
will probably lead to a better understanding of the link between (Euler-Mahonian) permutation and partition statistics.

\vspace{0.5cm}

\textbf{Acknowledgments.} The author would like to thank Professor D. Foata for suggesting 
to use partially commutative monoid in the study of partition statistics and 
 Einar Steingr\'\i msson for his help and encouragement during the preparation of this paper.



\end{document}